\documentclass[10pt]{amsart}

\usepackage{latexsym,amsfonts,amssymb,epsfig,verbatim}
\usepackage{amsmath, latexsym, graphics, textcomp,psfrag}
\usepackage{float}
\usepackage{amsthm}

\usepackage[all,2cell,dvips]{xy}

\newcommand{\nc}{\newcommand}
\nc{\dmo}{\DeclareMathOperator}
\nc{\nt}{\newtheorem}

\nc{\ds}{\displaystyle}
\nc{\ens}{\ensuremath}

\theoremstyle{plain}
\nt{thm}{Theorem}
\nt*{main}{Theorem}
\nt{lem}[thm]{Lemma}
\nt{cor}[thm]{Corollary}
\nt*{cor*}{Corollary}
\nt{prop}[thm]{Proposition}
\nt*{q}{Question}
\nt{fact}[thm]{Fact}
\nt*{claim}{Claim}
\nt{conjecture}[thm]{Conjecture}

\dmo{\mcg}{Mod}
\dmo{\sym}{Sym}

\nc{\mods}{\mcg(S)}
\nc{\Mod}{\mcg}

\nc{\Z}{\mathbb{Z}}
\nc{\R}{\mathbb{R}}
\nc{\Q}{\mathbb{Q}}
\nc{\B}{\mathcal{B}}
\nc{\K}{\mathcal{K}}
\renewcommand{\H}{\mathcal{H}}
\nc{\SE}{\textrm{SE}}
\dmo{\sy}{Sp}
\dmo{\SL}{SL}
\nc{\spz}{\sy(2g,\Z)}
\nc{\spnz}{\sy(2n,\Z)}
\nc{\spmz}{\sy(2m,\Z)}
\nc{\spgz}{\sy(2g,\Z)}
\nc{\slmz}{\SL(2m,\Z)}
\nc{\I}{\mathcal{I}}

\nc{\p}[1]{\medskip\paragraph{{\bf #1}}}
\nc{\margin}[1]{\marginpar{\scriptsize #1}}

\nc{\bpf}{\begin{proof}}
\nc{\epf}{\end{proof}}

\include{diagram}

\begin{document}

\input{epsf.sty}

\title{A homological recipe for pseudo-Anosovs}

\author{Dan Margalit}

\author{Steven Spallone}

\address{Dan Margalit \\ Department of Mathematics\\ University of
  Utah\\ 155 S 1440 East \\ Salt Lake City, UT 84112-0090.  Email: margalit@math.utah.edu}

\address{Steven Spallone \\ Department of Mathematics\\ Purdue University\\
150 N. University Street\\West Lafayette, IN 47907-2067.  Email: sspallon@math.purdue.edu}

\thanks{The first author gratefully acknowledges support from the National Science Foundation.}



\keywords{pseudo-Anosov, mapping class groups}

\subjclass[2000]{Primary: 20F36; Secondary: 57M07}

\maketitle


\begin{abstract}We give a simple explicit construction of
  pseudo-Anosov mapping classes using an improvement of the
  homological criterion of Casson--Bleiler.\end{abstract}

\section{Introduction}

Starting from a homological criterion of Casson--Bleiler, we give a
simple and explicit construction of pseudo-Anosov elements of the
mapping class group.  
In the exposition, we relegate much of the standard background material about
mapping class groups to other references; see, e.g., \cite{cb,fm}.

Let $S$ be a surface of genus $g$.  The \emph{mapping class group} of $S$, denoted $\Mod(S)$, is the group of isotopy classes of homeomorphisms of $S$.  The action of $\Mod(S)$ on $H_1(S,\Z)$ is symplectic, and in fact $\Mod(S)$ surjects onto the integral symplectic group; we call this map $\Psi$.  The kernel of $\Psi$ is called the \emph{Torelli group} for $S$, and is denoted $\I(S)$.  When $S$ has at most one boundary component, $H_1(S,\Z) \cong \Z^{2g}$, and we have the following short exact sequence.
\[ 1 \to \I(S) \to \mods \stackrel{\Psi}{\to} \spz \to 1 \]

Even given the Nielsen--Thurston classification, which says that every mapping class is periodic, reducible, or pseudo-Anosov, it is not at all obvious how to write down examples of pseudo-Anosov mapping classes.  Thurston gave a simple explicit construction, in particular giving the first examples of pseudo-Anosov elements of $\I(S)$ (Nielsen conjectured that no such mapping classes exist).

Since that time, various constructions of pseudo-Anosov mapping classes have been given in the literature; we do not attempt to give a survey here.
We do point out two advantages to our construction.  First, our examples are very simple to write down as a product of Dehn twists and ``handle switches'' (the latter can in turn be written as simple products of Dehn twists).  Second, since we are using a homological criterion, our construction is robust in the sense that if $f$ is one of our examples and $f' \in \I(S)$, then $ff'$ is pseudo-Anosov.

\p{Bounding triple maps and handle switches.}  Our construction involves two types of elementary mapping classes: bounding triple maps and handle switches.  Let $S$ be a surface of genus $g$ with at most one boundary component, and let $\{x_i,y_i\}$ be a standard symplectic basis for $H_1(S,\Z)$.  

For each $1 \leq i \leq g$, choose a mapping class $\B_i$ as follows.
\[
\B_i = 
\begin{cases}
T_{x_i+x_{i+1}}^{-1}T_{x_i}T_{x_{i+1}} & i < g, \ i \mbox{ odd} \\
T_{y_i+y_{i+1}}T_{y_i}^{-1}T_{y_{i+1}}^{-1}& i < g, \ i \mbox{ even} \\
T_{x_i}^{-1} & i = g, \ i \mbox{ odd} \\
T_{y_i} & i = g, \ i \mbox{ even} \\
\end{cases}
\]
As usual, $T_a$ denotes the (left) Dehn twist about a curve $a$ (we
are confusing homology classes with curves in the definition; any
curve in the given homology class will suffice). Each $\B_i$ with $i <
g$ is what we call a \emph{bounding triple map}, for each can be
realized as a product of Dehn twists about three mutually disjoint
curves, no pair of which separates $S$, but the union of which does.

By a \emph{handle switch} $\H_{i}$, we will mean any mapping class which achieves the following action on $H_1(S,\Z)$. 
\[
(x_i,y_i,x_{i+1},y_{i+1}) \mapsto
\begin{cases}
(-y_{i+1},x_{i+1},-y_i,x_i) & i < g, \ i \mbox{ odd} \\
(-y_{i+1},x_{i+1},-y_i,x_i) & i < g, \ i \mbox{ even} \\
\end{cases}\]\[
(x_i,y_i) \mapsto
\begin{cases}
(-y_i,x_i) & i = g, \ i \mbox{ odd} \\
(y_i,-x_i) & i = g, \ i \mbox{ even} \\
\end{cases}
\]
We give two ways of realizing a handle switch in
Section~\ref{section:mapping classes}.

\p{Symmetrization of polynomials.} Our construction requires one more
ingredient.  We define a
function $\sym:\Z[x] \to \Z[x]$ via the following formula.
\[ \sym(q)(x) = x^{\deg(q)} \cdot q(x+\frac{1}{x}) \]
The use of this function was suggested by Zrinka Despotovic.

\p{Statement of the construction.} We now are ready to describe a typical pseudo-Anosov mapping class yielded by our construction.  Let $S$, $\B_i$, and $\H_i$ be as above.

Let $p$ be any prime number, and let $a_0, \dots, a_{g-1}$ be integers
with $a_0$ not divisible by $p$ and $|a_{g-1}p|>2g$ (recall $g$ is
the genus of $S$).  Consider the polynomial
\[ q(x)=x^g+a_{g-1}px^{g-1} + \cdots + a_0p \]
and define integers $c_i$ via the following formula.
\[ \sym(q)(x) = x^{2g} + c_1x^{2g-1} + \cdots + c_{2g-1}x + 1 \]
Finally, let $t_i = (-1)^ic_i$, and consider the following mapping
class. 
\[ f = \B_1^{t_1}\H_1 \cdots \B_{g-1}^{t_{g-1}}\H_{g-1}\B_g^{t_g}\H_g \]

\begin{main}
Each element of the coset $f\I(S)$ is pseudo-Anosov.
\end{main}

\p{Questions.}  Building on a construction of Thurston, Penner gave
the following easy-to-state (and powerful) construction.  Let
$A=\{a_i\}$ and $B=\{b_i\}$ be two multicurves which fill a surface
$S$, and let $f$ be any mapping class obtained as a product of
positive powers of the $T_{a_i}$ and negative powers of the $T_{b_i}$
(each $a_i$ and $b_i$ must appear at least once).  Penner's theorem is
that $f$ is pseudo-Anosov.  Penner has conjectured that every
pseudo-Anosov mapping class has a power which is given by this
construction.  Thus, it is natural to ask whether or not this is true
for our examples.  Even if we write one such $f$ in ``Penner form,''
it is not clear how to write the other elements of $f\I(S)$ in this
way.

For any of our examples $f$, we would like to know if there is
an element of $f\I(S)$ whose stable and unstable foliations are
orientable.  In that case, the dilatation is exactly the spectral
radius of the matrix $\Psi(f)$.  It may be possible to
gain insight into the spectrum of dilatations of $\Mod(S)$ in this
way.  For instance, Farb--Leininger--Margalit showed that the minimal
dilatation pseudo-Anosov in $\Mod(S)$ does not lie in any of the level
$m$ congruence subgroups for $m \geq 3$ \cite{flm}; our construction
yields examples which do not lie in any such subgroup.  Chris
Leininger has pointed out the following two facts: only finitely many
$\Mod(S)$ conjugacy classes in a $\I(S)$ coset can have orientable
foliations (since in this case the dilatations are all the same, and it
is known that there are only finitely many conjugacy classes of
pseudo-Anosov elements with dilatation less than any given constant);
also, by finding a ``Casson--Bleiler polynomial'' with no real roots,
one can find examples of cosets where there are no representatives
with orientable foliations.

\p{Acknowledgements}  The authors are grateful to Ian Agol, Mladen Bestvina, Joan Birman, Zrinka Despotovic, Sam Evens, Benson Farb, Misha Kapovich, Chris Leininger, Lee Mosher, Dan Silver, and Susan Williams for their comments and help.


\section{The homological criterion}

In ``Automorphisms of surfaces after Nielsen and Thurston,'' Casson--Bleiler proved the following fact \cite[Lemma 5.1]{cb}.

\begin{prop}
\label{thm:cb hom crit}
Let $f \in \mods$, let $q_f(x)$ be the characteristic
polynomial for $\Psi(f)$.  If $q_f(x)$ is irreducible over $\Z$, has no roots
of unity as zeros, and is not a polynomial in $x^k$ for $k > 1$, then
$f$ is pseudo-Anosov.
\end{prop}

We remark that the only example to which Casson--Bleiler apply this homological criterion is a product of Dehn twists which (as they point out) is a case of the Penner construction mentioned earlier.
Below, we give a refinement of Proposition~\ref{thm:cb hom crit} which
lends itself better to an explicit construction of pseudo-Anosov
mapping classes; indeed, the examples given in the introduction do not
in general satisfy the hypotheses of Proposition~\ref{thm:cb hom crit}.  

We first need to introduce some terminology pertaining to
polynomials.  A (integer) \emph{symplectic polynomial} is any polynomial which
is the characteristic polynomial of some element of $\spnz$.  It is a
fact that an even degree integer polynomial
\[ q(x) = a_{2n} x^{2n} + \cdots + a_1 x + a_0 \]
is a symplectic polynomial if and only if it is monic ($a_{2n}=1$) and
palindromic ($a_i = a_{2n-i}$ for $0 \leq i \leq n$).  We say that a
symplectic polynomial is \emph{symplectically irreducible} if it is
not the product of two nontrivial symplectic polynomials.

Upon examination of the proof of Theorem~\ref{thm:cb hom crit} (see
\cite{cb}), we see that we can replace the word ``irreducible'' with
``symplectically irreducible''; this fact was pointed out to us by
Mladen Bestvina.  Also, as noticed by Silver--Williams, we can replace the assumption that $p_f(x)$ has no roots of unity as zeros with the assumption that $q_f(x)$ is cyclotomic \cite{sw}.

\begin{prop}[Homological criterion]
\label{thm:hom crit refined}
Let $S$ be a closed surface of genus at least 2, let $f \in \mods$,
and let $q_f(x)$ be the characteristic polynomial for $\Psi(f)$.  If
$q_f(x)$ is symplectically irreducible, is not a cyclotomic
polynomial, and is not a polynomial in $x^k$ for $k > 1$, then $f$ is
pseudo-Anosov.
\end{prop}

The hypotheses of Proposition~\ref{thm:hom crit refined} can be
thought of as hypotheses on the mapping class $f$, the matrix
$\Psi(f)$, or the polynomial $q_f(x)$.  We will say that a mapping
class, matrix, or polynomial \emph{satisfies the homological
criterion} if it (or the associated polynomial) satisfies the
hypotheses of the proposition.

Briefly, Proposition~\ref{thm:hom crit refined} is explained as
follows.  If a mapping class $f$ is periodic, then $\Psi(f)$ is
periodic, and so $q_f(x)$ has only roots of unity as zeros.  If $f$ is
reducible and $f^k$ fixes a nonseparating curve, then $q_f$ has a root
of unity as a zero. If $f$ is reducible and $f$ fixes a collection of
separating curves, then either $f$ preserves two distinct subsurfaces
of positive genus (in which case $q_f$ is symplectically reducible) or
$f$ transitively permutes a collection of subsurfaces which span
$H_1(S,\Z)$ (in which case the matrix has several identity blocks and
so $q_f$ is a polynomial in $x^k$ for some $k > 1$).  One also needs
the remark that a symplectically irreducible polynomial which has a
root of unity as a zero is necessarily cyclotomic.  Thus, if $q_f$
satisfies the homological criterion, then $f$ is neither periodic nor
reducible; by the Nielsen--Thurston classification, $f$ is
pseudo-Anosov.


\section{The construction}
\label{section:construction}

We will use the homological criterion to give a construction
of pseudo-Anosov mapping classes in the following way:

\begin{enumerate}
\item For every positive integer $g$, we give an explicit construction of
  infinitely many symplectic polynomials of degree $2g$, each of which
  satsifies the homological criterion.
\item Given any symplectic polynomial $q(x)$ of degree $2g$, we give
  an element $A_q$ of $\spz$ with characteristic polynomial $q(x)$.
\item For each such element $A_q$ of $\spz$, we give an infinite
  collection of mapping classes whose action on $H_1(S,\Z)$ is given
  by $A_q$.
\end{enumerate}

Each of the three steps is detailed below.  In
Section~\ref{section:genus 2}, we give a complete description of the
situation for genus 2 surfaces.

\subsection{Polynomials}

The first step is to find polynomials which satisfy the homological criterion. 
Recall the function $\sym: \Z[x] \to \Z[x]$ defined in the introduction.
We see that if $q(x)$ is monic, then $\sym(q)(x)$ is monic, of even degree, and is palindromic, i.e., $\sym(q)(x)$ is symplectic.

It is an elementary fact that the function
$\sym:\Z[x] \to \Z[x]$ is multiplicative, and injective.  Further, the monic elements surject onto the symplectic polynomials.  We thus have the following.

\begin{lem}
\label{fact:symp irred}
If $q(x) \in \Z[x]$ is monic, then $\sym(q)$ is symplectically
irreducible if and only if $q(x)$ is irreducible over $\Z$.
\end{lem}

From the definition of $\sym$, we see that $\zeta$ is a root of $\sym(q)$ if and only if $\zeta+\frac{1}{\zeta}$ is a root of $q$.  The following is a consequence.
 
\begin{lem}
\label{fact:root 2}
If $\sym(q)(x)$ is a cyclotomic
polynomial, then each root of $q(x)$ has absolute value bounded above by $2$.
\end{lem}

By a straightforward computation, we see that the coefficient of $x^{g-1}$ in a degree $g$ polynomial $q(x)$ is equal to the coefficients of $x$ and $x^{2g-1}$ in $\sym(q)(x)$.  We thus have the following.

\begin{lem}
\label{lem:simple}
If $q(x) \in \Z[x]$ is a monic degree $g$ polynomial, and the coefficient of
$x^{g-1}$ is nonzero, then $\sym(q)(x)$ is not a polynomial in $x^k$
for any $k > 1$.
\end{lem}

\begin{prop}
\label{prop:poly}
Let $q(x) \in \Z[x]$ be an irreducible polynomial of the following form.
\[ q(x) = x^g + a_{g-1}x^{g-1}+  \cdots + a_1 x + a_0 \]
If $|a_{g-1}|>2g$, then
$\sym(q)(x)$
satisfies the homological criterion.
\end{prop}

\bpf

By Lemma~\ref{fact:symp irred}, $\sym(q)(x)$ is symplectically
irreducible.  The condition on $a_{g-1}$ implies that some root of $q$
has absolute value greater than $2$ (indeed, $-a_{g-1}$ is the sum of
the roots of $q$); so by Lemma~\ref{fact:root 2} we see 
that $\sym(q)(x)$ is not a cyclotomic polynomial.
Finally, by Lemma~\ref{lem:simple}, $\sym(q)(x)$ is not a polynomial in
$x^k$ for any $k > 1$.
\epf

We can now use the Eisenstein criterion to construct explicit families of polynomials which satisfy the homological criterion.

\begin{prop}
Let $p$ be a prime number, and let $a_{g-1}, \dots, a_0$ be integers with
$|a_{g-1}p|>2g$ and $a_0$ not divisible by $p$.  If we set
\[ q(x)=x^g+a_{g-1}px^{g-1}+ \cdots + a_0p \]
then $\sym(q)$ satisfies the homological criterion.
\end{prop}

\subsection{Matrices}
\label{section:matrices}

We give the definition of the symplectic group and give two types of
basic elements: symplectic elementary matrices and symplectic
elementary permutation matrices.  We then use the theory of rational
canonical forms to give a symplectic matrix with any given symplectic
polynomial as its characteristic polynomial, and to factor the matrix
into a product of matrices of the above two types.

\p{Symplectic matrices.} To establish notation, we set
\[ J = \left(\begin{array}{ccccc}
             0 & 0 & \cdots & 0 & 1 \\
             0 & 0 & \cdots & -1 & 0 \\
             \vdots & \vdots & \ddots & \vdots & \vdots\\
             0 & 1 & \cdots & 0 & 0 \\
             -1 & 0 & \cdots & 0 & 0  \\
          \end{array} \right)    \]
and say that a matrix $A$ is symplectic when the following holds.
\[ A^T J A = J \]

Let $\{e_i\}$ be the standard basis vectors for $\Z^{2g}$.  Using the
inner product defined by $J$, we see that $\langle e_i, e_{2g-i+1}
\rangle = (-1)^{i+1}$ and all other inner products are zero; thus the $\{e_i\}$, in the order
$e_1,e_{2g},e_{2g-1},e_2,e_3,e_{2g-2},\dots$ give a symplectic basis
for $\Z^{2g}$.  We denote these vectors by $x_1,y_1,x_2,y_2,\dots,x_g,y_g$,
respectively.  Let $\sigma$ be the permutation of $1, \dots, 2g$
which is the product of the (commuting) transpositions $i
\leftrightarrow 2g-i+1$ for $1 \leq i \leq g$; note that
$\{e_i,e_{\sigma(i)}\} = \{x_i,y_i\}$.

There are two types of \emph{elementary symplectic
  matrices}.   If $i=\sigma(j)$, then the symplectic elementary matrix $\SE_{i,j}$ is given by the row
operation
\[ e_i \mapsto e_i + e_j \]
(that is, add the $j^{\mbox{\tiny{th}}}$ row to the
$i^{\mbox{\tiny{th}}}$ row; in general, when we write a row operation,
the unmentioned rows are fixed).  If $i \neq \sigma(j)$, then the elementary symplectic matrix $\SE_{i,j}$ is given
by the following row operations.
\begin{eqnarray*}
e_i  &\mapsto& e_i + e_j \\
e_{\sigma(j)}  & \mapsto & e_{\sigma(j)}  +(-1)^{i+j+1} e_{\sigma(i)}
\end{eqnarray*}

It is a fact that the $\SE_{i,j}$ generate the symplectic group $\spgz$.

Let $N_{i,j}$ denote the matrix given by the following row operation.
\[ (e_i,e_j) \mapsto (-e_j,e_i) \]
We then have the \emph{symplectic elementary permutation matrices} $W_{i,j}$.  If $i = \sigma(j)$, set
\[ W_{i,j} = N_{i,j} \]
and otherwise
\[ W_{i,j} = N_{i,j}N_{\sigma(j),\sigma(i)} \]

\p{Matrices with a given characteristic polynomial.}
Let $q(x)$ be the element $x^{2g}+c_{1}x^{2g-1} + \cdots + c_{2g}$ of
$\Z[x]$, and consider the following matrix.
\[ A_q =     \left(\begin{array}{rrrrr}
             -c_{1} & -c_{2} & \cdots & -c_{2g-1} & -c_{2g} \\
             1 & 0 & \cdots & 0 & 0 \\
             0 & 1 & \cdots & 0 & 0 \\
             \vdots & \vdots & \ddots & \vdots & \vdots\\
             0 & 0 & \cdots & 1 & 0  \\
          \end{array} \right)    \]

\begin{lem}
The matrix $A_q$ has characteristic
polynomial $q(x)$.
\end{lem}

If the determinant of $A_q$ (namely $c_{2g}$) is equal to $1$, one can
obtain $A_q$ as a product of elementary matrices and permutation
matrices in a simple way.  Let $t_i=(-1)^{i}c_i$.  Write $X_i$ for
the elementary matrix corresponding to the row operation $e_i \mapsto
e_i + t_i e_{i+1}$, $N_i$ for $N_{i,i+1}$, and $Y_i$ for $X_iN_i$.
We have the following theorem \cite[\S 7.4, \S 7.8]{S}.

\begin{thm}
\label{rcf}
Using the notation above, we have
\[ A_q=Y_1 Y_2 \cdots Y_{2g-1} \]
Moreover, if the $Y_i$ are multiplied in a different order, the result
is conjugate to $A_q$, and thus has the same characteristic polynomial (namely, $q(x)$).
\end{thm}

Consider the following expression.
\[ B_q= (Y_1Y_{2g-1})(Y_2Y_{2g-2}) \cdots (Y_{g-1}Y_{g+1}) Y_g \]
Note that
\[ Y_iY_{2g-i}=X_iN_iX_{2g-i}N_{2g-i}=(X_iX_{2g-i})(N_iN_{2g-i}) \]
since $N_i$ and $X_{2g-i}$ commute.  If $q(x)$ is symplectic, then the
matrices $X_iX_{2g-i}$, $N_iN_{2g-i}$, and $Y_g=X_gN_g$ are symplectic.
Indeed, for $i < g$, we have
\begin{eqnarray*}
X_{i}X_{2g-i} &=& \SE_{i,i+1}^{t_i} \\
N_iN_{2g-i} &=& W_{i,i+1}\\
\end{eqnarray*}
and also
\begin{eqnarray*}
X_g &=& \SE_{g,g+1}^{t_g} \\
N_g &=& W_{g,g+1}
\end{eqnarray*}
Thus we have expressed $B_q$, which is conjugate to $A_q$ (Theorem~\ref{rcf}), explicitly
as a product of symplectic elementary matrices and symplectic elementary
permutation matrices.
\[ B_q = \SE_{1,2}^{t_1} W_{1,2}  \SE_{2,3}^{t_2} W_{2,3} \cdots \SE_{g-1,g}^{t_{g-1}} W_{g-1,g} \SE_{g,g+1}^{t_g} W_{g,g+1} \]
For our purposes, this is convenient, because
it is not difficult to find simple, explicit mapping classes whose actions on homology
are given by $\SE_{i,i+1}$ and $W_{i,i+1}$.

\p{Remark} One could also apply the ``Steinberg section'' [4] to the symplectic
group to immediately obtain the same expression for $B_q$.



\subsection{Mapping classes}
\label{section:mapping classes}

\begin{figure}[htb]
\psfrag{x1}{$x_1$}
\psfrag{x2}{$x_2$}
\psfrag{xg}{$x_g$}
\psfrag{y1}{$y_1$}
\psfrag{y2}{$y_2$}
\psfrag{yg}{$y_g$}
\centerline{\includegraphics[scale=.5]{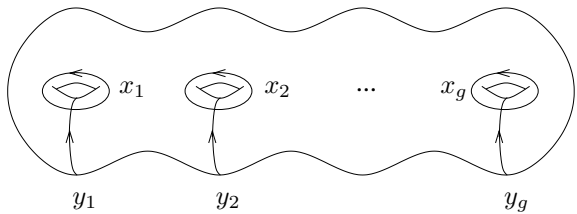}}
\caption{The standard symplectic basis for $H_1(S,\Z)$.}
\label{figure:basis}
\end{figure}

To this point, we have given explicit examples of symplectic
polynomials which satisfy the homological criterion, and we have given
explicit matrices with the given polynomials as characteristic
polynomials.  Moreover, we have given these matrices as products of the
symplectic elementary matrices $\SE_{i,j}$ and symplectic elementary permutation matrices
$W_{i,j}$.  Thus, to complete our construction of pseudo-Anosov
mapping classes, it suffices to exhibit mapping classes whose actions
on $H_1(S,\Z)$ are given by the $\SE_{i,j}$ and $W_{i,j}$.  We will use the standard basis for $H_1(S,\Z)$ shown in Figure~\ref{figure:basis}.

If $x$ is an element of $H_1(S,\Z)$, then we write $T_x$ for the \emph{Dehn twist} (to the left) about any simple closed curve representing $x$.  As a result of this convention, we will not be specifying particular elements of $\Mod(S)$ below, but rather cosets of $\I(S)$, from which particular elements are easily extracted.

For $1 \leq i,j \leq 2g$, set $t = (-1)^{j+1}$, and define
\[
\B_{i,j} =
\begin{cases}
T_{e_i}^t & i=\sigma(j) \\
(T_{e_i+e_{\sigma(j)}}T_{e_i}^{-1}T_{e_{\sigma(j)}}^{-1})^t & i \neq \sigma(j)
\end{cases}
\]

(Recall that the $e_i$ are the standard basis elements for $\Z^{2g} \cong H_1(S,\Z)$ and are identified with the $x_i$ and $y_i$ as in Section~\ref{section:matrices}.)  In the case $j=i+1$, we get the mapping classes $\B_i$ given in the introduction.

\begin{lem}
\label{fact:se as maps}
For any $i$ and $j$, $\Psi(\B_{i,j}) = \SE_{i,j}$.
\end{lem}

By a change of coordinates, it suffices to check Lemma~\ref{fact:se as
maps} for a genus 2 surface and the 4 possibilities for $i$ and $j$
even and odd.

We remark that Lemma~\ref{fact:se as maps}, combined with the fact that the $\SE_{i,j}$ generate $\spgz$, constitutes a proof of the classical theorem that $\Psi:\mcg(S) \to \spgz$ is surjective.

For $i=\sigma(j)$, $W_{i,j}$ corresponds to the following signed row permutation.
\[ (e_i,e_{\sigma(i)}) \mapsto (-e_{\sigma(i)},e_i) \]
Such a mapping class can be visualized as a quarter-turn of a genus
one subsurface spanned homologically by $e_i$ and $e_{\sigma(i)}$.  As
a product of Dehn twists, this is the following mapping class.
\[ (T_{e_i}T_{e_{\sigma(i)}}T_{e_i})^{(-1)^{i+1}} \]  When $i=g$, we
denote this map by $\H_{g}$ (to be consistent with the introduction).

If $i \neq \sigma(j)$, the matrix $W_{i,j}$ is given by the following signed row permutation.
\[ (e_i,e_{\sigma(i)},e_j,e_{\sigma(j)}) \mapsto (-e_j,e_{\sigma(j)},e_i,-e_{\sigma(i)}) \]
In order to construct a mapping class whose action on
homology is a permutation matrix, one can choose standard representatives for the $e_k$, cut
$S$ along each $e_k$ in order
to obtain a surface $S'$ with $g$ ``square'' boundary components, take
any homeomorphism of $S'$ where the action on the squares is given by
the action on homology, and then ``glue up'' to get a homeomorphism of
$S$.

If one prefers, it is also possible to find an explicit product of
Dehn twists which does the same job.  For simplicity, we consider the
case where $i$ and $j$ have different parity and neither is equal to
$g$ (this suffices for the construction).  Let $\H_i$ be the mapping
class
\[ ((T_{e_{\sigma(i)}}T_{e_i}T_{e_{\sigma(i)}+e_{i+1}}T_{e_{\sigma(i+1)}}T_{e_{i+1}})^3
(T_{e_{\sigma(i)}}T_{e_i})^3)^t\]  where $t=(-1)^{i+1}$ 
(the first
factor switches the two handles and the second factor turns one of the
handles by $180^\textrm{\tiny{o}}$).

\begin{lem}
\label{fact:weyls as maps}
For $1 \leq i \leq g$, $\Psi(\H_i) = W_{i,i+1}$.
\end{lem}

Again, this lemma can be checked by reducing the case of $4 \times 4$
matrices (genus 2 surfaces).

By Lemmas~\ref{fact:se as maps} and~\ref{fact:weyls as maps}, one can verify that the construction
given in the introduction is correct.  The most general statement
we will make is as follows.

\begin{main}
Let $S$ be a surface of genus $g$ with at most one boundary
component.  Let $q(x) = x^g + a_{g-1}x^{g-1} + \cdots + a_0$ be a
monic degree $g$ irreducible integral polynomial with $|a_{g-1}| >
2g$.  Define integers $c_i$ via
\[ \sym(q)(x) = x^{2g} + c_1x^{2g-1} + \cdots + c_{2g-1}x + 1 \]
and let $t_i = (-1)^ic_i$.  Let $f$ a
product of the elements of $\{(\B_i^{t_i}\H_i)\}_{i=1}^g$, where 
$\H_i$ and $\B_i$ are as in the introduction, and the product is taken
in any order.  In this case, each element of the coset $f\I(S)$ is
pseudo-Anosov.
\end{main}

\p{Remark.} Kanenobu has found, for each symplectic polynomial, a fibered link whose associated monodromy has that polynomial as the characteristic polynomial of the action of the monodromy on the fiber surface \cite{tk}.  However, we do not know how to explicitly understand these monodromies topologically on the surface.


\section{Genus 2}
\label{section:genus 2}

Let $S$ be the closed orientable surface of genus 2.  We proceed in
the same general outline as in Section~\ref{section:construction}.
That is, we first say exactly when a degree 4 symplectic polynomial
satisfies the homological criterion, then give an explicit product of
matrices for each symplectic polynomial, and then give an explicit and
simple product of Dehn twists corresponding to the product of matrices.

\begin{prop}
A degree 4 symplectic polynomial 
\[ q(x) = x^4 + a x^3 + b x^2 + a x +1 \]
satisfies the homological criterion
if and only if $a \neq 0$,
$(a,b)$ is not $(1,1)$ or $(-1,1)$, and
$a^2-4b+8$ is not a perfect square.
\end{prop}

\bpf

It is straightforward to check that $q(x) = \sym(h(x))$, where
\[ h(x) = x^2 + ax + (b-2) \]
By Lemma~\ref{fact:symp irred}, $q(x)$ is symplectically reducible if
and only if
$h(x)$ is reducible over $\Z$.  By the quadratic formula, this is
exactly when $-a \pm \sqrt{a^2-4b+8}$ is an even integer.  Note that this
is always the case when the discriminant is a perfect square.  By checking
the list of degree 4 cyclotomic polynomials, $q(x)$ is cyclotomic
exactly
when (a,b) is $(1,1)$, $(0,0)$, or $(-1,1)$.  Finally, $q(x)$ is a
polynomial in $x^4$ or $x^2$ exactly when $a = 0$.
\epf

As in Section~\ref{section:matrices}, the matrix
\[ A_q=
\left(\begin{array}{rrrr}
-a & -b & -a & -1 \\
1 & 0 & 0 & 0 \\
0 & 1 & 0 & 0 \\
0 & 0 & 1 & 0  \\
          \end{array} \right)
\]
has characteristic polynomial $q(x) = x^4 + a x^3 + b x^2 + a x +1$ (this can be checked by direct computation).

As per Section~\ref{section:matrices}, we consider the product
\[ (X_1X_3)(N_1N_3)(X_2)(N_2) = \SE_{1,2}^{-a}W_{1,3}\SE_{2,3}^bW_{2,3}
\]
which is equal to the following.
\[
\left(\begin{array}{rrrr}
1 & -a & 0 & 0 \\
0 & 1 & 0 & 0 \\
0 & 0 & 1 & -a \\
0 & 0 & 0 & 1  \\
          \end{array} \right)
\left(\begin{array}{rrrr}
0 & -1 & 0 & 0 \\
1 & 0 & 0 & 0 \\
0 & 0 & 0 & -1  \\
0 & 0 & 1 & 0 \\
          \end{array} \right)
\left(\begin{array}{rrrr}
1 & 0 & 0 & 0 \\
0 & 1 & b & 0 \\
0 & 0 & 1 & 0 \\
0 & 0 & 0 & 1  \\
          \end{array} \right)
\left(\begin{array}{rrrr}
1 & 0 & 0 & 0 \\
0 & 0 & -1 & 0 \\
0 & 1 & 0 & 0 \\
0 & 0 & 0 & 1  \\
          \end{array} \right)
\]
This last product is conjugate to $A_q$ and hence has characteristic
polynomial $q(x)$.  We are free to change the signs of the middle two
columns in the second and fourth matrices since this operation has order
2 and commutes with the third matrix.
\[
\left(\begin{array}{rrrr}
1 & -a & 0 & 0 \\
0 & 1 & 0 & 0 \\
0 & 0 & 1 & -a \\
0 & 0 & 0 & 1  \\
          \end{array} \right)
\left(\begin{array}{rrrr}
0 & 1 & 0 & 0 \\
1 & 0 & 0 & 0 \\
0 & 0 & 0 & -1  \\
0 & 0 & -1 & 0 \\
          \end{array} \right)
\left(\begin{array}{rrrr}
1 & 0 & 0 & 0 \\
0 & 1 & b & 0 \\
0 & 0 & 1 & 0 \\
0 & 0 & 0 & 1  \\
          \end{array} \right)
\left(\begin{array}{rrrr}
1 & 0 & 0 & 0 \\
0 & 0 & 1 & 0 \\
0 & -1 & 0 & 0 \\
0 & 0 & 0 & 1  \\
          \end{array} \right)
\]

These four matrices, from left to right, are the images of the mapping
classes $(T_{x_1+x_2}T_{x_1}^{-1}T_{x_2}^{-1})^a$,
$(T_{y_2}T_{x_2}T_{y_1+y_2}T_{x_1}T_{y_1})^3$, $T_{y_2}^b$, and
$T_{x_2}T_{y_2}T_{x_2}$, respectively.

\begin{prop}
If $a$ and $b$ are integers with the properties
that $a \neq 0$, $(a,b)$ is not $(1,1)$ or $(1,-1)$ and $a^2-4b+8$
is not a perfect square, and $f$ is any element of $\I(S)$, then
the mapping class
\[ (T_{x_1+x_2}T_{x_1}^{-1}T_{x_2}^{-1})^a(T_{y_2}T_{x_2}T_{y_1+y_2}T_{x_1}T_{y_1})^3(T_{y_2})^b(T_{x_2}T_{y_2}T_{x_2}) f \]
is pseudo-Anosov.
\end{prop}

\bibliographystyle{plain}
\bibliography{homcrit}

\end{document}